\documentclass[12pt]{article}
\usepackage[latin1]{inputenc}
\usepackage{amsfonts}
\usepackage{amsmath}
\usepackage{amssymb}
\usepackage{amsfonts}
\usepackage{amscd}
\usepackage{bbm}
\usepackage{graphicx}
\usepackage{theorem}
\theoremstyle{break}
\newtheorem{The}{Theorem}
\newtheorem{Lem}{Lemma}
\newtheorem{Def}{Definition}
\newtheorem{Pro}{Proposition}
\title{Selmer's multiplicative algorithm}
\author{By\\ J. Christopher Kops}
\date{}
\begin{document}
\maketitle
\begin{abstract}
The behavior of the multiplicative acceleration of Selmer's algorithm is widely unknown and no general result on convergence has been detected yet. Solely for its $2$-dimensional, periodic expansions exist some results on convergence and approximation due to Fritz Schweiger $\cite{Sch2005}$. In this paper we show that periodic expansions of any dimension do in fact converge and that the coordinates of the limit points are rational functions of the largest eigenvalue of the periodicity matrix.
\end{abstract}
\section{Introduction}
There are at least two ways to approach a theory of multidimensional continued fractions (MCFs). One is more of a geometric nature, while another concentrates on multidimensional continued fractions which can be described by a set of $(n+1)\times(n+1)$-matrices. The letter set includes amongst others the Jacobi-Perron algorithm, as well as algorithms of Brun and Selmer. Each of these generalizes the matrices $\begin{pmatrix} a&1 \\ 1&0\end{pmatrix}$ associated with continued fractions to higher dimensions in order to achieve an equivalent of Lagrange's theorem, but for cubic or higher roots. However it still remains a challenge, whether one of these algorithms applied to an arbitrary $n$-tuple $x=(\alpha_1,\dots,\alpha_n)$ whereas $1,\alpha_1,\dots,\alpha_n$ are linearly independent over $\mathbb{Q}$ and belong to a real number field of degree $n+1$ eventually becomes periodic. \\
Despite this lack of knowledge a lot of discoveries have been made concerning periodicity and approximation characteristics of multidimensional continued fractions. For further information the books of Brentjes $\cite{Bre1981}$ and Schweiger $\cite{Sch2000}$, both named "multidimensional continued fractions" are highly recommended.
\subsection{Fibred systems}\label{sub:fib}
In this subsection we first introduce the notion of fibred systems and illustrate its characteristics in the context of continued fractions. Thus we are able to define multidimensional continued fractions by a set of matrices on such a fibred system. %Beforehand we show that the one-dimensional continued fraction algorithm is indeed a fibred system.\\
\begin{Def}[fibred system]
Let $B$ be a set and $T\colon B\to B$ be a map. The pair $(B,T)$ is called a fibred system if the following conditions are satisfied:
\begin{enumerate}
\item There is a finite or countable set I (called the digit set).
\item There is a map $k\colon B\to I$. Then the sets
\begin{equation*}
B(i)=k^{-1}\{i\}=\{x\in B: k(x)=i\}
\end{equation*}
form a partition of $B$, hence $\bigcup_{i\in I}B(i)=B$.
\item The restriction of $T$ to any $B(i)$ is an injective map.
\end{enumerate}
(The partition $\{B(i): i\in I\}$ is called the time-1-partition)%MCF S.1
\end{Def}
\begin{bfseries}
Example
\end{bfseries}
We consider the continued fraction algorithm. As we know to each irrational $x\in[0,1]$ there corresponds a continued fraction expansion
\begin{equation*}
x=[a_0,a_1,a_2,\dots]=a_0+\cfrac1{a_1 + \cfrac1{a_2+\dots}},
\end{equation*}
where we restrict $a_1,a_2,\dots$ to be positive integers. Now we set $B:=[0,1[$ and obtain for $x\neq 0$ the map 
\begin{equation*}\begin{aligned}
T\colon B&\to B\\
 x&\mapsto \frac{1}{x}-a(x), \quad a(x):=\lfloor{x^{-1}}\rfloor.
\end{aligned}\end{equation*}
For $x=0$ however we arrive at a stop and thus get $Tx=0$.
Further by setting $I:=\mathbb{N}_0$ we obtain the map $i\colon B\to I$ defined by 
\begin{equation*}
i\colon x\mapsto a(x).
\end{equation*}
Then the set $I:=\mathbb{N}_0$ is countable, while the sets $B(i)$ form a partition of B. Now we set $y:=Tx=\frac{1}{x}-a(x)$ and restrict $T$ to one of the sets $B(i), i\in I$. Then for $x_1,x_2\in B(i)$ with $Tx_1= Tx_2$ we immediately obtain $x_1=x_2$, since $a(x_1)=a(x_2)=i$. Hence the pair $(B,T)$ indicates a fibred system.\\\\
With the knowledge of fibred systems we are now able to define multidimensional continued fractions on the aforementioned systems simply by a set of matrices.
\begin{Def}[multidimensional continued fraction]\label{Def:mcf}%MCF S.2
The fibred system $(B,T)$ is called (multidimensional) continued fraction if
\begin{enumerate}
\item $B$ is a subset of ${\mathbb{R}}^n$.
\item For every digit $k\in I$ there is an invertible matrix $\alpha=\alpha(k)=((A_{ij}))$, \\ $0\leq i, j \leq n$, such that $y=Tx,\,x\in B(k)$, is given as
\begin{equation*}
y_i=\frac{A_{i0}+\sum_{j=1}^n{A_{ij}x_j}}{A_{00}+\sum_{j=1}^n{A_{0j}x_j}}.
\end{equation*}
\end{enumerate}
\end{Def}
In particular, we are interested in the inverse matrix of such a multidimensional continued fraction, as hereby we obtain an expression of $x$ via its image under the map $T$.
\begin{Def}\label{eq: def6}%MCF S.4
If $(B,T)$ is a (multidimensional) continued fraction, then we denote the inverse matrix of $\alpha(k)$ by ${\beta(k)}=((B_{ij}))$, $0\leq i, j\leq n$. Then we define for $s\geq 1$:
\begin{equation*}
\beta(k_1,\dots,k_s):=\beta(k_1)\dots\beta(k_s)=((B^{(s)}_{ij})), \quad0\leq i, j\leq n,
\end{equation*}
whereas $B^{(1)}_{ij}=B_{ij}$.
Then $y=T^sx$ is equivalent to
\begin{equation*}
x_i=\frac{B^{(s)}_{i0}+\sum^n_{j=1}B^{(s)}_{ij}y_j}{B^{(s)}_{00}+\sum^n_{j=1}B^{(s)}_{0j}y_j},\quad1\leq i \leq n.
\end{equation*}
\end{Def}
\subsection{The concept of cylinders}\label{sub:cyl} 
In the previous subsection we saw that the digit set $I$ of a fibred system causes a partition of the set $B$ into the subsets $B(i),i\in I$, that is
\begin{equation*}
\bigcup_{i\in I}B(i)=B.
\end{equation*}
Simultaneously we didn't demand that $T$ restricted to $B(i),i\in I$ is a surjective map. Therefore a partition of $TB(i),i\in I$ could be of interest as well and we eventually arrived at the concept of cylinders.
\begin{Def}[cylinder]\label{Def:cylinder}%Ergodic Theory S.3+MCF S.1
The cylinder of rank $s$, defined by the digits $i_1,\dots,i_s$ is the set
\begin{equation*}
B(i_1,\dots, i_s):=B(i_1)\cap T^{-1}B(i_2)\cap\ldots\cap T^{-(s-1)}B(i_s)=\{x: i_1(x)=i_1,\dots, i_s(x)=i_s\}.
\end{equation*}
Such a cylinder $B(i_1,\dots,i_s)$ is called proper (or full) if $T^sB(i_1,\dots,i_s)=B$. 
\end{Def}
\begin{Pro}\label{Pro:proper}
All cylinders of arbitrary rank $s,\, s\in \mathbb{N}$ are full if all cylinders of rank $1$ are full.
\end{Pro}
\begin{bfseries}
Beweis:
\end{bfseries}
Assume that all cylinders of rank $1$ are full, then we get 
\begin{equation*}
TB(i)=B
\end{equation*}
for all $i\in I$. The rest follows from the Definition \eqref{Def:cylinder} of cylinders, since for $(i_1,\dots,i_s)\in I^s$ we obtain
\begin{equation*}
TB(i_1,\dots,i_s)=B(i_2,\dots,i_s).
\end{equation*}
$\hfill\fbox{}$
\clearpage
\section{Selmer's algorithm}
In connection with Brun's algorithm, Selmer published 1961 a variation of its subtractive version called the subtractive algorithm of Selmer (SSA). But instead of subtracting the second biggest initial value from the biggest like Brun did, he chose to subtract the smallest initial value from the biggest. This may at first seem to be a marginal deviation, but it actually implicates a fundamental change. As at some point in the expansion one inevitably ends up in the absorbing set $D:=B(n-1)\cup B(n)$. 
\subsection{Subtractive version}
Let
\begin{equation*}
\Delta^{n+1}:=\{b=(b_0,b_1,\dots,b_n):b_0\geq b_1\geq\dots\geq b_n\geq0\},
\end{equation*}
then we define
\begin{equation*}
\sigma b:=(b_0-b_n,b_1,\dots,b_n).
\end{equation*}
There is an index $i=i(b),\, 0\leq i\leq n$, such that
\begin{equation*}
\pi\sigma b:=(b_1,b_2,\dots,b_i,b_0-b_n,\dots,b_n)\in\Delta^{n+1}.
\end{equation*}
Further let
\begin{equation*}
B^n:=\{(x_1,\dots,x_n):\, 1\geq x_1\geq\dots\geq x_n\geq 0\}.
\end{equation*}
Then with the help of the projection $p:\Delta^{n+1}\rightarrow B^n$ defined by 
\begin{equation*}
p(b_0,b_1,\dots,b_n)=\left(\frac{b_1}{b_0},\dots,\frac{b_n}{b_0}\right)
\end{equation*}
we obtain the map $T:B^n\rightarrow B^n$ which makes the diagramm
\[ \begin{CD} 
\Delta^{n+1} @>\pi\sigma>> \Delta^{n+1} \\ 
@V{p}VV @VV{p}V \\ 
B^n @>>T> B^n 
\end{CD} \] 
commutative.
\subsection{Fibred system and absorbing set}
In this subsection we prove that the SSA intrinsically represents a fibred system. Afterwards we shortly point to the absorbing set of Selmer's algorithm. As a matter of fact, the restriction of the SSA to this absorbing set coincides with an algorithm of M\"onkemeyer. Hence for further references we point to $\cite{Mon1954}$, $\cite{Pan2008}$ and $\cite{Sch2000}$.
\begin{Pro}\label{Pro:partition}
The Partition
\begin{equation*}
B(j):=\left\{x\in B^n:i(p^{-1}x)=j\right\},\quad 0\leq j\leq n,
\end{equation*}
makes $(B^n,T)$ a fibred system.
\end{Pro}
\begin{bfseries}
Proof:
\end{bfseries}
We set $y:=Tx$ and calculate the map $T$ resp. the inverse branches $V(j)$, $0\leq j\leq n$, as follows: \\
\fbox{$j=0$}
\begin{align*}
y_k&=\frac{x_k}{1-x_n},\quad1\leq k\leq n,\quad& x_k&=\frac{y_k}{1+y_n},\quad 1\leq k\leq n.
\intertext{$\fbox{$1\leq j\leq n-1$}$}
y_k&=\frac{x_{k+1}}{x_1},\quad1\leq k\leq j-1,\quad& x_1&=\frac{1}{y_j+y_n},\\
y_j&=\frac{1-x_{n}}{x_1},\quad \quad& x_k&=\frac{y_{k-1}}{y_j+y_n},\quad2\leq k\leq j,\\
y_k&=\frac{x_{k}}{x_1},\quad j+1\leq k\leq n,\quad& x_k&=\frac{y_{k}}{y_j+y_n},\quad j+1\leq k\leq n.
\intertext{$\fbox{$j=n$}$}
y_k&=\frac{x_{k+1}}{x_1},\quad1\leq k\leq n-1,&\quad x_1&=\frac{1}{y_{n-1}+y_n},\\
y_n&=\frac{1-x_n}{x_1},&\quad x_k&=\frac{y_{k-1}}{y_{n-1}+y_n},\quad2\leq k\leq n.
\end{align*}
Therefore we obtain:
\begin{enumerate}
\item The digit set $I:=\{0,\dots,n\}$ is a countable set.
\item The map $i\colon B^n\to I$ defined by $i\colon x\mapsto i(p^{-1}x)$ causes a partition of $B^n$, since evidently $\bigcup_{j\in I}B(j)=B^n$.
\item If we restrict $T$ to any $B(j)$, then we obtain for $y^{'}=y^{''}$ and $0\leq j\leq n$:
\begin{enumerate}
\item[$\fbox{$j=0$}$] As $y_n^{'}=y_n^{''}$ it follows that $x_n^{'}=x_n^{''}$ and thus $x^{'}=x^{''}$.\\
\item[$\fbox{$1\leq j\leq n-1$}$] As $y_n^{'}=y_n^{''}$ and $y_j^{'}=y_j^{''}$ we obtain $x_1^{'}=x_1^{''}$ and hence $x^{'}=x^{''}$.\\
\item[$\fbox{$j=n$}$] As $y_n^{'}=y_n^{''}$ and $y_{n-1}^{'}=y_{n-1}^{''}$ we get $x_1^{'}=x_1^{''}$ and therefore $x^{'}=x^{''}$.
\end{enumerate}
Consequently the restriction of $T$ to any $B(j)$ is an injective map.
\end{enumerate}
$\hfill\fbox{}$
\begin{The}[absorbing set]%[MCF - S. 55]
\label{Sat:absorb}
Let $D:=\{x\in B^n: x_{n-1}+x_n\geq 1\}$. Then $D$ is an absorbing set, i.e.
\begin{enumerate}
\item $TD=D$
\item For almost every $x\in B^n$ there is an $N=N(x)$, such that $T^Nx\in D$.
\end{enumerate}
\end{The}
\begin{bfseries}
Proof:
\end{bfseries}
A proof of this theorem can be found in $\cite[p. 55]{Sch2000}$ and verifies that
\begin{equation*}
D=B(n-1)\cap B(n).
\end{equation*}
$\hfill\fbox{}$
\subsection{A case of periodicity}
As the question of periodicity is by far the most interesting one, we now give an example of a periodic SSA.
\begin{Def}[periodic continued fraction]\label{Def:periodic}%MCF S.111 
The multidimensional continued fraction of $x$ is called periodic if there are numbers $m\geq 0$, $p\geq 1$ such that $T^{m+p}x=T^mx$.
\end{Def}
\begin{bfseries}
Example (SSA)
\end{bfseries}
We consider the tuple $x:=(x_1,x_2)=\left(\sqrt[3]{4}-1, \sqrt[3]{2}-1\right)$ and apply Selmer's algorithm. Then we obtain
\begin{align*}
Tx&=\left(\frac{\sqrt[3]{4}-1}{2-\sqrt[3]{2}}, \frac{\sqrt[3]{2}-1}{2-\sqrt[3]{2}}\right)\\
T^2x&=\left(\frac{3-2\sqrt[3]{2}}{\sqrt[3]{4}-1}, \frac{\sqrt[3]{2}-1}{\sqrt[3]{4}-1}\right)\\
&\vdots\\
T^{30}x&=\left(\frac{54-29\sqrt[3]{2}-11\sqrt[3]{4}}{30\sqrt[3]{4}+13\sqrt[3]{2}-64}, \frac{24\sqrt[3]{2}+3\sqrt[3]{4}-35}{30\sqrt[3]{4}+13\sqrt[3]{2}-64}\right)\\
T^{31}x&=\left(\frac{27\sqrt[3]{4}-11\sqrt[3]{2}-29}{54-29\sqrt[3]{2}-11\sqrt[3]{4}}, \frac{24\sqrt[3]{2}+3\sqrt[3]{4}-35}{54-29\sqrt[3]{2}-11\sqrt[3]{4}}\right).
\end{align*}
and an easy calculation shows that $T^{31}x=Tx$. Hence the SSA for $x:=(x_1,x_2)=\left(\sqrt[3]{4}-1, \sqrt[3]{2}-1\right)$ becomes periodic with a preperiod of length 1 and a period with length 30.
%\clearpage
\section{Multiplicative algorithms}
As for Brun's algorithm the multiplicative version really causes an acceleration of expansions, the multiplicative version of Selmer's algorithm in general doesn't. In addition none of its cylinders are full. Thus it's more adequate to denote it as a mere division algorithm. Under these circumstances it seems unlikely that the algorithm provides convergent expansions or approximations that are competitive in the field of multidimensional continued fractions. However, it does.
\subsection{Selmer's division algorithm}
Let $\Delta^{n+1}:=\{b=(b_0,b_1,\dots,b_n):\, b_0\geq b_1\geq\dots\geq b_n\geq0\}$. Then we define
\begin{equation*}
\delta b:=(b_0-kb_n,b_1,\dots,b_n),\, k:=\left[\frac{b_0}{b_n}\right].
\end{equation*}
Since $b_n\geq b_0-kb_n$ we get $\pi\colon \Delta^{n+1}\to\Delta^{n+1}$ defined by
\begin{equation*}
\pi\delta b:=(b_1,\dots,b_n,b_0-kb_n).
\end{equation*}
Now let $B^n:=\{(x_1,\dots,x_n):\, 1\geq x_1\geq\dots\geq x_n\geq 0\}$. With the help of the projection $p\colon \Delta^{n+1}\to B^{n}$ defined by 
\begin{equation*}
p(b_0,b_1,\dots,b_n)=\left(\frac{b_1}{b_0},\dots,\frac{b_n}{b_0}\right)
\end{equation*}
we finally get the bottom map
\begin{align*}
S&(x_1,\dots,x_n)=\left(\frac{x_2}{x_1},\dots,\frac{x_n}{x_1},\frac{1-kx_n}{x_1}\right),
\end{align*}
which makes the diagramm
\[ \begin{CD} 
\Delta^{n+1} @>\pi\sigma>> \Delta^{n+1} \\ 
@V{p}VV @VV{p}V \\ 
B^n @>>S> B^n 
\end{CD} \] 
commutative.\\
Hence the multiplicative version of Selmer's algorithm (MSA) is given by
\begin{align*}
S&\colon B^n\to B^n\\
S&(x_1,\dots,x_n)=\left(\frac{x_2}{x_1},\dots,\frac{x_n}{x_1},\frac{1-kx_n}{x_1}\right).
\end{align*}
\subsection{The fibred system}
In this subsection we simply prove that the MSA represents a fibred system.
\begin{Pro}
The partition
\begin{equation*}
B(k):=\left\{x\in B^n:\frac{1}{k+1}<x_n\leq\frac{1}{k}\right\},\quad k=1,2,\dots
\end{equation*}
makes $(B^n,S)$ a fibred system.
\end{Pro}
\begin{bfseries}
Proof:
\end{bfseries}
If $x_n=0$, we simply restrict to $B^{n-1}:=\{(x_1,\dots,x_{n-1}):\, 1\geq x_1\geq\dots\geq x_{n-1}\geq 0\}$. Then from $k:=[x_n^{-1}]$ we immediately get $k\in\mathbb{N}$ and thus:
\begin{enumerate}
\item The digit set $I:=\mathbb{N}$ is a countable set.
\item The map $i\colon B^n\to I$ amounts with $i\colon x\mapsto k:=[x_n^{-1}]$ to a partition of $B^n$, since evidently $\bigcup_{k\in I}B(k)=B^n$.
\item If we restrict $T$ to $B(k)$, then for $y{'}=y{''}$ and thus $y_n{'}=y_n{''},\,y_{n-1}{'}=y_{n-1}{''}$ we obtain $x_1{'}=x_1{''}$ and hence $x{'}=x{''}$. Therefore the restriction of $T$ to $B(k)$ is an injective map for all $k\in\mathbb{N}$.
\end{enumerate}
$\hfill\fbox{}$
\subsection{Cylinders and time-1-partition}
\begin{figure}[h]
\centering
\includegraphics{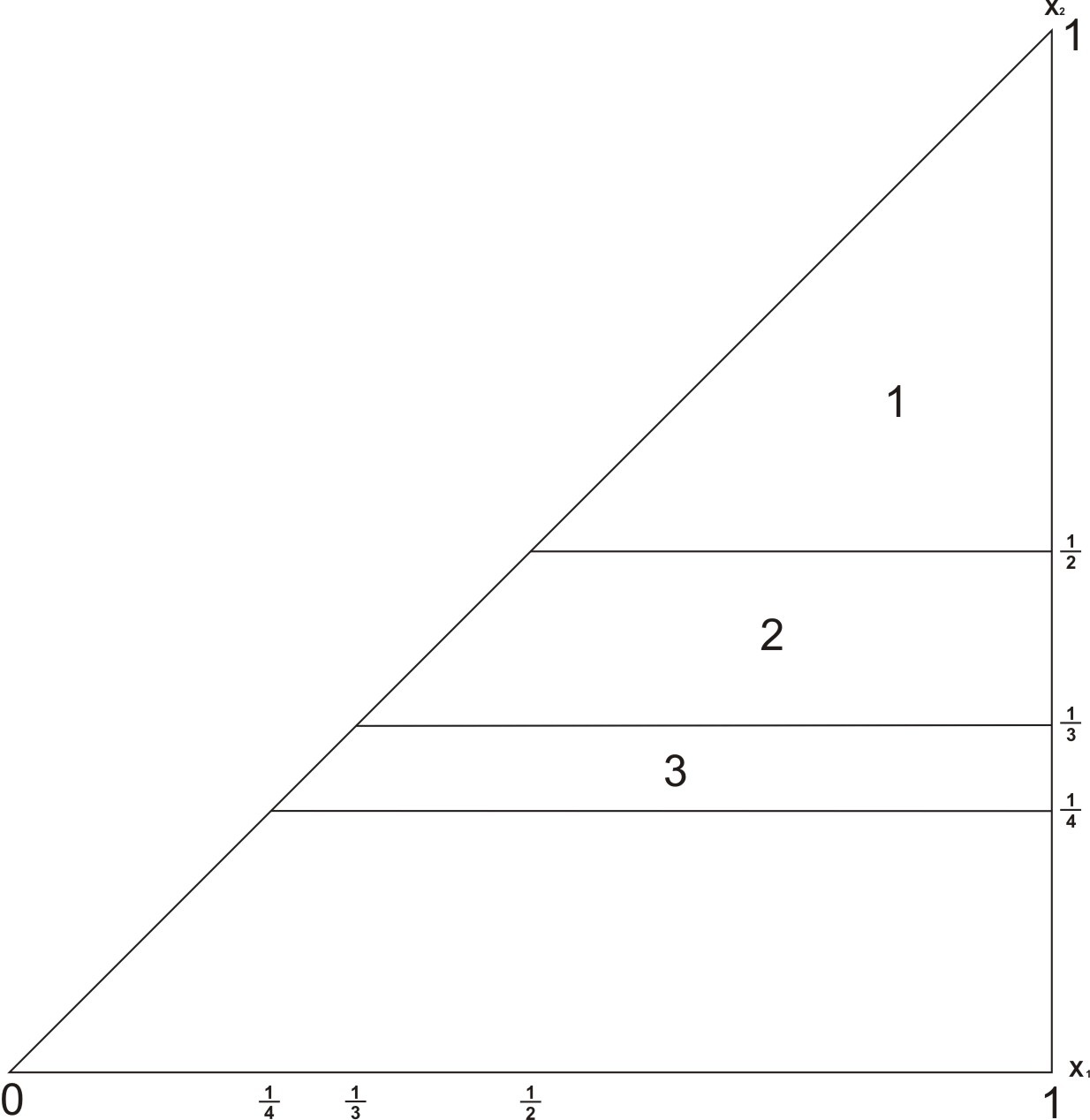} 
\caption{The time-1-partition of the set $B^2$ by the 2-dimensional MSA, where $k\in\mathbb{N}$ indicates the associated cylinder $B(k)$.}\label{fig:fig1}
\end{figure}

Since $k=[x_n^{-1}]$ the pair $(B^n,S)$ is a fibred system with cells
\begin{equation*}
B(k):=\left\{x\in B^n:\frac{1}{k+1}<x_n\leq\frac{1}{k}\right\},\quad k=1,2,\dots.
\end{equation*}
These $B(k)$ denote cylinders of rank 1 and indicate in our case convex sets with vertices
\begin{align*}
&\left(1,\dots,1,\frac{1}{k}\right),\left(1,\dots,1,\frac{1}{k+1}\right)\\
&\left(1,\dots,\frac{1}{k},\frac{1}{k}\right),\left(1,\dots,\frac{1}{k+1},\frac{1}{k+1}\right)\\
&\quad\quad\vdots\\
&\left(\frac{1}{k},\dots,\frac{1}{k}\right),\left(\frac{1}{k+1},\dots,\frac{1}{k+1}\right)
\end{align*}
only depending on $k\in\mathbb{N}$. By the way we obtain $S(\frac{1}{k+1},\dots,\frac{1}{k+1})=(1,\dots,1)$ for all cylinders $B(k),\,k\in\mathbb{N}$.\\
In order to achieve more clearness we restrict our attention to the $2$-dimensional case. Thus we now consider the set 
\begin{equation*}
B^2=\{1\geq x_1\geq x_2\geq 0\}.
\end{equation*}
and since $1\geq x_1\geq x_2$ und $\frac{1}{k}\geq x_2>\frac{1}{k+1}$ we obtain the convex set $B(k)$ with vertices
\begin{equation*}
\left(1,\frac{1}{k}\right),\left(1,\frac{1}{k+1}\right),\left(\frac{1}{k},\frac{1}{k}\right),\left(\frac{1}{k+1},\frac{1}{k+1}\right).
\end{equation*}
These cells $B(k)$ form a partition of the set $B^2$ that can be easily illustrated by means of figure \eqref{fig:fig1}.
\\However none of the cylinders $B(k)$ are full, as they are mapped under $S$ onto the convex set with vertices
\begin{equation*}
\left(\frac{1}{k},0\right),\left(\frac{1}{k+1},\frac{1}{k+1}\right),\left(1,0\right),\left(1,1\right).
\end{equation*}
Hence $SB(k)\subset SB(k+1)$ for all $k\in\mathbb{N}$ and additionally $SB(k)$ is not a union of cylinders of rank $1$.
\subsection{The matrices}
Notice that if we set $y:=Sx$, then
\begin{align*}
x_1&=\frac{1}{ky_{n-1}+y_n}\\
x_i&=\frac{y_{i-1}}{ky_{n-1}+y_n},\quad 2\leq i\leq n.
\end{align*}
Thus, according to the definition of a multidimensional continued fraction, the associated $(n+1)\times(n+1)-$matrices of the MSA are given by
\begin{equation*}\beta(k):=
\begin{pmatrix} 0 & \dots  & k & 1 \\
1 & 0 & \dots & 0\\
\vdots & \ddots & \ddots& \vdots\\
0 & \dots & 1 & 0
\end{pmatrix}.
\end{equation*}
and in dimension $n=3$ by
\begin{equation*}
\begin{pmatrix} 
0 & 0 & k & 1 \\
1 & 0 & 0 & 0\\
0 & 1 & 0 & 0\\
0 & 0 & 1 & 0
\end{pmatrix}.
\end{equation*}
We see immediately, that expanding the determinant along the last column leads to
\begin{equation}\label{eq:det}
\det\beta(k)=(-1)^{n}\det\mathbbm{1}=\pm 1.
\end{equation}
Now we define the matrices $\beta^{(s)}(k_1,\dots,k_s)$ in common notation as
\begin{equation*}
\beta^{(s)}(k_1,\dots,k_s):=\beta(k_1)\dots\beta(k_s)=
\begin{pmatrix}
B_0^{(s-n+1)} & \dots & B_0^{(s-1)} & B_0^{(s)} & B_0^{(s-n)} \\
\vdots&\vdots &  & & \vdots\\
B_n^{(s-n+1)} & \dots & B_n^{(s-1)} & B_n^{(s)} & B_n^{(s-n)} 
\end{pmatrix}.
\end{equation*}
Hence for $s\geq 0$ we obtain the relation
\begin{equation*}
B_i^{(s+1)}=k_{s+1}B_i^{(s-n+1)}+B_i^{(s-n)},\quad i=0,\dots,n,
\end{equation*}
whereas $\beta^{(0)}$ denotes the unit matrix. If we set $y=S^sx$ and $k_i=k(S^{i-1})$, $1\leq i\leq s$, then we find that
\begin{equation}\label{eq:relation}
x_i=\frac{B_i^{(s-n+1)}+y_1B_i^{(s-n+2)}+\dots+y_{n-1}B_i^{(s)}+y_nB_i^{(s-n)}}{B_0^{(s-n+1)}+y_1B_0^{(s-n+2)}+\dots+y_{n-1}B_0^{(s)}+y_nB_0^{(s-n)}},\quad i=1,\dots,n.
\end{equation}
\subsection{Periodic expansions for Selmer's division algorithm}
In this subsection we eventually prove the convergence of the periodic MSA and quote a simple example of periodicity. Since this proof allows us to apply a variety of other theorems to Selmer's multiplicative algorithm, we subsequently mention a few of them. However for completeness and further information we refer to \cite{Sch2000}.
\subsubsection{Weak convergence of the periodic MSA}
\begin{Def}[weak convergence]
The multidimensional continued fraction is weakly convergent if for every $x\in B$
\begin{align*}
\lim_{s\rightarrow\infty}{\left(\frac{B^{(s)}_{10}}{B^{(s)}_{00}},\dots,\frac{B^{(s)}_{n0}}{B^{(s)}_{00}}\right)}&=x.
\intertext{If for arbitrary $g$ with $0\leq g\leq n$ we get}
\lim_{s\rightarrow\infty}{\left(\frac{B^{(s)}_{1g}}{B^{(s)}_{0g}},\dots,\frac{B^{(s)}_{ng}}{B^{(s)}_{0g}}\right)}&=x,\quad 
\end{align*}
for every $x\in B$ we call the multidimensional continued fraction uniformly weakly convergent.
\end{Def}
\begin{The}\label{Sat:konv}
Assume that the algorithm of $x=(x_1,\dots,x_n)$ eventually becomes periodic with period length $p$. Then
\begin{equation*}
\lim_{s\to\infty}\left(\frac{B_1^{(s)}}{B_0^{(s)}},\dots,\frac{B_n^{(s)}}{B_0^{(s)}}\right)=x.
\end{equation*}
Thus the periodic, multiplicative algorithm of Selmer is weakly convergent and even uniformly weakly convergent. %\cite[p. 6]{Sch2000}
\end{The}
\begin{bfseries}
Proof:
\end{bfseries}
Clearly we can assume that the expansion is purely periodic with period length $p$.\\
Let
\begin{equation*}
M:=\beta^{(p)}(k_1,\dots,k_p)=
\begin{pmatrix}
B_0^{(p-n+1)} & \dots & B_0^{(p-1)} & B_0^{(p)} & B_0^{(p-n)} \\
\vdots & \vdots & & &\\
B_n^{(p-n+1)} & \dots & B_n^{(p-1)} & B_n^{(p)} & B_n^{(p-n)} 
\end{pmatrix}.
\end{equation*}
Note that
\begin{equation*}
M^k=
\begin{pmatrix}
B_0^{(kp-n+1)} & \dots & B_0^{(kp-1)} & B_0^{(kp)} & B_0^{(kp-n)} \\
\vdots & \vdots & & &\\
B_n^{(kp-n+1)} & \dots & B_n^{(kp-1)} & B_n^{(kp)} & B_n^{(kp-n)}
\end{pmatrix}.
\end{equation*}
Then the characteristic polynomial of $M$ is given as
\begin{equation*}
\chi_M(t):=\det{(t\mathbbm{1}-M)}=t^{n+1}-b_nt^n-\dots-b_1t-b_0,
\end{equation*}
whereas $b_0=(-1)^{n-1}\det{(M)}$. Furthermore we denote its eigenvalues by $\rho_0,\rho_1,\dots,\rho_n$.
Assume that the matrix $M$ has entries equal to $0$, then with the help of Lemma \eqref{Lem:positiv} we get $m(n)$, such that $M^m$ is a positive matrix for all $m\geq m(n)$ and we could continue with the period length $mp$. Thus we can assume, that $M$ is a positive matrix and the period length remains $p$.\\
Now the Perron-Frobenius theorem \cite[Seite 112]{Sch2000} allows us to arrange the eigenvalues in the way
\begin{equation*}
\rho_0>|\rho_1|\geq\dots\geq|\rho_n|,\quad \rho_0>1,
\end{equation*}
whereas $\rho_0$ is a simple and positive root of $\chi_M(t)$. Further with the help of the famous Caley-Hamilton Theorem we obtain
\begin{equation*}
M^{n+1}-b_nM^n-\dots-b_1M-b_0=0
\end{equation*}
and by multiplying this with $M^k\beta(k_1,\dots,k_j)$ we see that for the entries of the matrices the relations
\begin{equation*}
B_i^{((k+n+1)p+j)}-b_nB_i^{((k+n)p+j)}-\dots-b_1B_i^{((k+1)p+j)}-b_0B_i^{(kp+j)}=0, 
\end{equation*}
hold for $0\leq i\leq n$ and $0\leq j< p$. Now we apply Theorem 41 of \cite[Seite 114]{Sch2000} and therefore obtain the general solution
\begin{equation}\begin{aligned}\label{eq:goldenboy}
B_i^{(kp+j)}=\quad&d(i,j)\rho_0^k\\
&+b_{10}(i,j)\rho_1^k+b_{11}(i,j)\begin{pmatrix} k \\ 1 \end{pmatrix}\rho_1^{k-1}+\dots +b_{1m_1}(i,j)\begin{pmatrix} k \\ m_1-1 \end{pmatrix}\rho_1^{k-m_1+1}\\
&+\dots\\
&+b_{s0}(i,j)\rho_s^k+b_{s1}(i,j)\begin{pmatrix} k \\ 1 \end{pmatrix}\rho_s^{k-1}+\dots +b_{sm_s}(i,j)\begin{pmatrix} k \\ m_s-1 \end{pmatrix}\rho_s^{k-m_s+1},
\end{aligned}\end{equation}
whereas $m_1,\dots,m_s$ are the multiplicities of the roots $\rho_1,\dots,\rho_s$ and thus $1+m_1+\dots+m_s=n+1$. If the start values $B_i^{(kp+j)}$ are given for all $k=0,\dots,n$, then the solution sequence 
$(B_i^{(kp+j)}), k\geq 1$ is uniquely determined.\\
Now we will consider the terms in equation \eqref{eq:goldenboy} more precisely. Clearly from $\rho_0>|\rho_1|\geq\dots\geq|\rho_n|$ and $\rho_0>1$ follows
\begin{equation*}
\rho_0^k>\left|\rho_{\gamma}^{k-\mu}\right|,\quad 1\leq \gamma\leq s,\quad 0\leq \mu\leq m_{\gamma}-1.
\end{equation*}
Furthermore for all $\gamma,\,1\leq\gamma\leq s$ and $\varepsilon>0$ the relations
\begin{align*}
\begin{pmatrix} k+1 \\ \mu \end{pmatrix}&=\begin{pmatrix} k \\ \mu \end{pmatrix}\frac{k+1}{k+1-\mu},\quad 1\leq\mu\leq m_{\gamma}-1\\
(1+\varepsilon)^{k+1}&=(1+\varepsilon)^{k}(1+\varepsilon)
\end{align*}
hold for all $k\in\mathbb{N}$ and quite evidently follows
\begin{align*}
\frac{k+1}{k+1-\mu}&\xrightarrow{k\to\infty} 1,\quad 1\leq\mu\leq m_{\gamma}-1\\
(1+\varepsilon)&>1. 
\end{align*}
Hence there exists a $k^{'}$, such that for all $k\geq k^{'}(\varepsilon)$ the inequality $\frac{k+1}{k+1-\mu}<(1+\varepsilon)$ holds and consequently with increasing $k,\,k\geq k^{'}(\varepsilon)$ the term $(1+\varepsilon)^k$ grows faster than $\begin{pmatrix}k\\\mu\end{pmatrix}$. Thus there exists a $k(\varepsilon)\geq k^{'}(\varepsilon)$, such that for all $k\geq k(\varepsilon)$ the inequality
\begin{equation*}
\begin{pmatrix}k\\\mu\end{pmatrix}<(1+\varepsilon)^k\quad\forall\text{ $\mu$ with } 1\leq\mu\leq m_{\gamma}-1.
\end{equation*}
is valid. Since $\rho_0>1$, there is a $\varepsilon>0$, such that $\rho_0=1+\varepsilon$. Consequently there is also a
$k^*(\varepsilon)>k(\varepsilon)$, such that for all $k\geq k^*(\varepsilon)$ the inequality
\begin{equation*}
\rho_0^k>\begin{pmatrix}k\\\mu\end{pmatrix}\left|\rho_{\gamma}^{k-\mu}\right|,\quad 0\leq\mu\leq m_{\gamma}-1.
\end{equation*}
holds. So we derive from equation \eqref{eq:goldenboy} for $d(i,j)\neq 0$ the limit
\begin{equation}\label{eq:ice}
\lim_{k\to\infty}\frac{B_i^{(kp+j)}}{\rho_0^k}=d(i,j),
\end{equation}
as $n,\,n\in\mathbbm{N}$ is a finite integer and $b_{\gamma \mu}$ represents a constant for all $\gamma ,\mu$ with $1\leq \gamma\leq s,0\leq \mu\leq m_{\gamma}$, respectively.\\
Now we consider the recursion relations
\begin{equation*}\begin{aligned}
B_i^{(s+n)}&=k_{s+n}B_i^{(s)}+B_i^{(s-1)}\\
B_i^{(s+n+1)}&=k_{s+n+1}B_i^{(s+1)}+B_i^{(s)}.
\end{aligned}\end{equation*}
It follows, that if $B_i^{(s)}\gg \rho_0^k$ then $B_i^{(s+j)}\gg \rho_0^k$ for all $j\in\mathbb{N}$, $j\geq n\cdot(n-1)$. Since $d(i,j)$ doesn't depend on $k$ and $n$ is finite, it follows that if $d(i,j)\neq 0$ for some $j$, then $d(i,j)\neq 0$ for all $j$, $0\leq j<p$.\\
We know that for the trace of $M^k$ the relation
\begin{equation}\label{eq:beliebig}
B_0^{(kp-n+1)}+\dots+B_{n-1}^{(kp)}+B_n^{(kp-n)}=\rho_0^k+\dots+\rho_n^k
\end{equation}
holds and since this equation \eqref{eq:beliebig} is valid for arbitrarily big $k$ and $n$ is finite, there is at least one summand $B_0^{(kp-n+1)},\dots,B_{n-1}^{(kp)},B_n^{(kp-n)}$, for whose related $d(i,j)$ follows $d(i,j)\neq 0$ by equation \eqref{eq:goldenboy}. Hence for this $i$ we obtain $d(i,j)\neq 0$ for all $j,0\leq j< p$.\\
Now we consider the relation for $x_i$ in equation \eqref{eq:relation}. As $y=T^sx\in B(k_{s+1}):=\{x\in B^n:\frac{1}{k_{s+1}+1}<x_n\leq\frac{1}{k_{s+1}}\}$ we get
\begin{equation*}
\frac{1}{k_{s+1}+1}\leq\frac{B_i^{(s-n+1)}+\dots+B_i^{(s)}+B_i^{(s-n)}}{B_0^{(s-n+1)}+\dots+B_0^{(s)}+B_0^{(s-n)}}\leq k_{s+1}+1,\quad i=1,\dots,n
\end{equation*}
and by defining $k^{*}:=\max(k_1,\dots,k_p)+1$ we eventually obtain 
\begin{equation}\label{eq:sumup}
\frac{1}{k^{*}}\leq\frac{B_i^{(s-n+1)}+\dots+B_i^{(s)}+B_i^{(s-n)}}{B_0^{(s-n+1)}+\dots+B_0^{(s)}+B_0^{(s-n)}}\leq k^{*},\quad i=1,\dots,n.
\end{equation}
Since equation \eqref{eq:sumup} is bounded as well for arbitrarily big $s$, but $d(i,j)\neq 0$ is valid for $0\leq j< p$ and a fixed $i\in\{0,\dots,n\}$, we obtain $d(i,j)\neq 0$ for all $0\leq i\leq n,\, 0\leq j<p$.\\
The equation
\begin{equation*}
M^k\begin{pmatrix}1 \\ x_1 \\ \vdots \\ x_n\end{pmatrix}=\lambda^k\begin{pmatrix}1 \\ x_1 \\ \vdots \\ x_n\end{pmatrix}
\end{equation*}
shows that for a positive eigenvalue $\lambda$ of the matrix $M$ the relation
\begin{equation}\label{eq:lambda}
B_0^{(kp-n+1)}+x_1B_0^{(kp-n+2)}+\dots+x_{n-1}B_0^{(kp)}+x_nB_0^{(kp-n)}=\lambda^k.
\end{equation}
holds.
But since equation \eqref{eq:lambda} is valid for arbitrarily big $k$ and $d(0,j)\neq 0$ for all $0\leq j<p$, we get $\lambda=\rho_0$. Thus by equation \eqref{eq:relation} we derive for $T^px=x$ the relation 
\begin{equation}\label{eq:revelation}
x_i=\rho_0^{-k}\left(B_i^{(kp-n+1)}+x_1B_i^{(kp-n+2)}+\dots+x_{n-1}B_i^{(kp)}+x_nB_i^{(kp-n)}\right)\quad 1\leq i\leq n.
\end{equation}
As this equation \eqref{eq:revelation} holds for arbitrarily big $k$ and $d(i,j)$ exists for all  $0\leq i\leq n$, $0\leq j< p$, the limit
\begin{equation}
x_i=\lim_{k\to\infty}\rho_0^{-k}\left(B_i^{(kp-n+1)}+x_1B_i^{(kp-n+2)}+\dots+x_{n-1}B_i^{(kp)}+x_nB_i^{(kp-n)}\right)
\end{equation}
exists for all $i,\,1\leq i\leq n$. With the help of equation \eqref{eq:ice} we obtain for all $i,\,1\leq i\leq n$ the relation
\begin{equation}\label{eq:solution}
x_i=d(i,p-n+1)+x_1d(i,p-n+2)+\dots+x_{n-1}d(i,0)+x_nd(i,p-n).
\end{equation}
Now we consider the equation $\beta^{2kp+j}=M^k\beta^{kp+j}$ and see that for the entries of the matrices the relations
\begin{equation*}
B_i^{2kp+j}=B_i^{(kp-n+1)}B_0^{(kp+j)}+\dots+B_i^{(kp)}B_{n-1}^{(kp+j)}+B_i^{(kp-n)}B_{n}^{(kp+j)}
\end{equation*}
hold for $1\leq i\leq n$ and $0\leq j\leq p-1$. Further we derive $\rho_0^k>|\rho_1^k|\geq\dots\geq|\rho_n^k|$ from $\rho_0>|\rho_1|\geq\dots\geq|\rho_n|$ and hence
\begin{equation}\label{eq:deees}
d(i,j)=d(i,p-n+1)d(0,j)+\dots+d(i,0)d(n-1,j)+d(i,p-n)d(n,j).
\end{equation}
follows.\\
Now we set 
\begin{equation}\label{eq:Salzburg}
x_i=\frac{d(i,j)}{d(0,j)},\quad 1\leq i\leq n
\end{equation}
in equation \eqref{eq:solution} and obtain equation \eqref{eq:deees} as a result. By equation 
\eqref{eq:det} we know, that the determinant of $M^k$ is unequal to zero. Thus equation \eqref{eq:solution} represents a system of $n$ equations in $n$ variables and $x_i$ is uniquely determined by equation \eqref{eq:Salzburg}. Since this is valid for all $j,\,0\leq j<p$ and $d(i,j)$ doesn't depend on $k$ we get
\begin{equation*}
\lim_{s\to\infty}\left(\frac{B_1^{(s)}}{B_0^{(s)}},\dots,\frac{B_n^{(s)}}{B_0^{(s)}}\right)=x.
\end{equation*}
$\hfill\fbox{}$
\begin{Lem}\label{Lem:positiv}
Let $M$ be a $(n+1)\times(n+1)-$matrix defined by 
\begin{equation*}M:=
\begin{pmatrix} 0 & \dots  & k & 1 \\
1 & 0 & \dots & 0\\
\vdots & \ddots & \ddots& \vdots\\
0 & \dots & 1 & 0
\end{pmatrix},
\end{equation*}
whereas $k\in\mathbb{N}$ represents a positive integer. Then there is a $p(n)\in\mathbb{N}$, such that $M^p$ is a positive matrix for all $p\geq p(n)$.
\end{Lem}
\begin{bfseries}
Proof:
\end{bfseries}
Notice that $M=E'+K$, whereas
\begin{equation*}E':=
\begin{pmatrix} 0 & \dots  & 0 & 1 \\
1 & 0 & \dots & 0\\
\vdots & \ddots & \ddots& \vdots\\
0 & \dots & 1 & 0
\end{pmatrix},\quad K:=
\begin{pmatrix} 0 & \dots  & k & 0 \\
0 & 0 & \dots & 0\\
\vdots & \ddots & \ddots& \vdots\\
0 & \dots & 0 & 0
\end{pmatrix}.
\end{equation*}
Furthermore $(E')^{n+1}=\mathbbm{1}$ represents the unit matrix and $K^2$ the zero matrix. Then apparently the relation
\begin{equation*}
K(E')^{n-1}=
\begin{pmatrix}
k &\cdots & 0
\\
\vdots & & \vdots
\\
0 & \cdots & 0
\end{pmatrix}
\end{equation*}
holds.\\
Now we set $p=n^2+1$ and consider $M^p$. Then we get $M^p=K^{*}+A$, whereas $K^{*}$ is given by
\begin{equation*}
K^{*}=
\begin{pmatrix}
k^{n-1} & \cdots & k & 1 & k^n \\
k^{n} & \ddots &  &  k & 1 \\
1 & \ddots & & & k \\
\vdots & \ddots & & & \vdots\\
k^{n-2} & \cdots & 1 & k^n & k^{n-1}
\end{pmatrix}
\end{equation*}
and $A\geq 0$. The shape of $K^{*}$ results from
\begin{align*}
K^{*}=(E')^{n^2+1}&+\sum_{j=0}^{n-1}\sum_{i=0}^{n}(E')^{i}[K(E')^{n-1}]^jK(E')^{n^2-jn-i}\\
=(E')^{n^2+1}&+\sum_{i=0}^{n}(E')^{i}K(E')^{n^2-i}\\
&+\sum_{i=0}^{n}(E')^{i}K(E')^{n-1}K(E')^{n^2-n-i}\\
&+\cdots\\
&+\sum_{i=0}^n(E')^{i}\underbrace{K(E')^{n-1} \dots K(E')^{n-1}}_{[K(E')^{n-1}]^{n-1}}K(E')^{n-i}.
\end{align*}
Hence $M^p$ has only posivite entries for all $p\geq n^2+1$ and therefore $\beta^{p}$ is a positive matrix, since $k\geq 1$.
$\hfill\fbox{}$
\subsubsection{Example of a periodic MSA}
Unfortunately the MSA of $x:=(\sqrt[3]{4},\sqrt[3]{2})$ doesn't become periodic within the first 40 steps of expansion.
Although we are unaware, whether or not periodicity eventually occurs in this expansion, there certainly are periodic expansions, some even of period length 1.\\
\begin{bfseries}
Example (MSA)
\end{bfseries}
We consider the tuple $x:=(x_1,x_2)=(\frac{\sqrt{5}-1}{2},\frac{3-\sqrt{5}}{2})$ and apply Selmer's multiplicative algorithm. Note that $x_2=x_1^2$. Then we obtain
\begin{align*}
Tx&=\left(\frac{3-\sqrt{5}}{\sqrt{5}-1}, \frac{-4+2\sqrt{5}}{\sqrt{5}-1}\right)
\intertext{and by multiplying each fraction with $\frac{\sqrt{5}+1}{\sqrt{5}+1}$ we obtain}
Tx&=\left(\frac{\sqrt{5}-1}{2},\frac{3-\sqrt{5}}{2}\right).
\end{align*}
Hence the MSA for $x:=(x_1,x_2)=(\frac{\sqrt{5}-1}{2},\frac{3-\sqrt{5}}{2})$ becomes periodic with a period of length 1.
\subsubsection{Some general results on the periodic MCF}
Since we proved in Theorem \eqref{Sat:konv} that the MSA is uniformly weakly convergent we can show that the coordinates of the limit points are rational functions of the largest eigenvalue of the periodicity matrix. Secondary we can apply a theorem on approximation properties of MCFs which can be traced back to Perron \cite{Per1907}. Their general forms for multidimensional continued fractions can be found in \cite[p. 115-119, 157-162]{Sch2000}.
\begin{The}
Assume that the algorithm of $x=(x_1,\dots,x_n)$ eventually becomes periodic with period length $p$. Then 
$x_1,\dots,x_n$ are rational functions in $\rho_0$, whereas $\rho_0$ denotes the largest eigenvalue of the characteristic polynomial of the periodicity matrix $\beta^{(p)}$. Therefore $x_1,\dots,x_n$ belong to a number field of degree $\leq n+1$.
\end{The}
\begin{bfseries}
Proof:
\end{bfseries}
Clearly we can assume that the expansion is purely periodic with period length $p$. Due to Theorem \eqref{Sat:konv} we know, that
\begin{equation*}
M^k\begin{pmatrix}1 \\ x_1 \\ \vdots \\ x_n\end{pmatrix}=\rho_0^k\begin{pmatrix}1 \\ x_1 \\ \vdots \\ x_n\end{pmatrix}
\end{equation*}
and thus get
\begin{equation*}
B_0^{(kp-n+1)}+x_1B_0^{(kp-n+2)}+\dots+x_{n-1}B_0^{(kp)}+x_nB_0^{(kp-n)}=\rho_0^k.
\end{equation*}
Therefore we can calculate $x_1,\dots,x_n$ as rational functions in $\rho_0$ from the equation
\begin{equation*}
M\begin{pmatrix}1 \\ x_1 \\ \vdots \\ x_n\end{pmatrix}=\rho_0\begin{pmatrix}1 \\ x_1 \\ \vdots \\ x_n\end{pmatrix}
\end{equation*}
$\hfill\fbox{}$\\\\
Below we imply by writing $a\ll b$ with $a,b\in\mathbb{R}$ that there is a constant $C\in\mathbb{N}$, such that $a\leq Cb$ is valid and accordingly by writing $a\gg b$ with $a,b\in\mathbb{R}$ that there is a constant $C\in\mathbb{N}$, such that the inequality $a\geq Cb$ holds.
\begin{The}
For $n$-dimensional, periodic algorithms we obtain for arbitrary $\varepsilon>0$ the estimates
\begin{equation*}
|B_0^{(pg)}x_i-B_i^{(pg)}|\ll |\rho_1(1+\varepsilon)|^g,\quad i=1,\dots,n.
\end{equation*}
But there is at least one pair $(i,j)$, such that the inequality
\begin{equation*}
|B_0^{(pg+j)}x_i-B_i^{(pg+j)}|\gg|\rho_1|^g
\end{equation*}
holds for infinitely many values of $g$.
\end{The}
\begin{bfseries}
Proof:
\end{bfseries}
Due to Theorem \eqref{Sat:konv} we get for all $i,\,1\leq i\leq n$ the relation
\begin{equation*}\begin{aligned}
B_0^{(pg)}x_i-B_i^{(pg)}=\quad&k(i,0)\rho_0^g\\
&+m_{10}(i,0)\rho_1^g+m_{11}(i,0)\begin{pmatrix} g \\ 1 \end{pmatrix}\rho_1^{g-1}+\dots +m_{1m_1}(i,0)\begin{pmatrix} g \\ m_1-1 \end{pmatrix}\rho_1^{g-m_1+1}\\
&+\dots\\
&+m_{s0}(i,0)\rho_s^g+m_{s1}(i,0)\begin{pmatrix} g \\ 1 \end{pmatrix}\rho_s^{g-1}+\dots +m_{sm_s}(i,0)\begin{pmatrix} g \\ m_s-1 \end{pmatrix}\rho_s^{g-m_s+1}
\end{aligned}\end{equation*}
and with the help of equation \eqref{eq:Salzburg} it follows, that $k(i,0)=0$ for all $i,\,1\leq i\leq n$. Thus the relation reduces to
\begin{equation*}\begin{aligned}
B_0^{(pg)}x_i-B_i^{(pg)}=\quad&m_{10}(i,0)\rho_1^g+m_{11}(i,0)\begin{pmatrix} g \\ 1 \end{pmatrix}\rho_1^{g-1}+\dots +m_{1m_1}(i,0)\begin{pmatrix} g \\ m_1-1 \end{pmatrix}\rho_1^{g-m_1+1}\\
+&\dots\\
+&m_{s0}(i,0)\rho_s^g+m_{s1}(i,0)\begin{pmatrix} g \\ 1 \end{pmatrix}\rho_s^{g-1}+\dots +m_{sm_s}(i,0)\begin{pmatrix} g \\ m_s-1 \end{pmatrix}\rho_s^{g-m_s+1}.
\end{aligned}\end{equation*}
By Theorem \eqref{Sat:konv} we know, that for any $\varepsilon>0$ there exists a $g(\varepsilon)$, such that for all $g\geq g(\varepsilon)$ the inequality
\begin{equation*}
\begin{pmatrix}g \\ \mu\end{pmatrix}<(1+\varepsilon)^g,\quad 1\leq\mu\leq m_{\gamma}-1,\quad1\leq \gamma\leq s.
\end{equation*}
holds. Since for the eigenvalues $|\rho_1|\geq\dots\geq|\rho_s|$ holds, $|\rho_1|^g\geq\dots\geq|\rho_s|^g$ is valid as well, where $s\leq n$. Hence there is a constant $c(x,\varepsilon)$, such that for all $g\geq g(\varepsilon, c)$ the inequality
\begin{equation*}
|B_0^{(pg)}x_i-B_i^{(pg)}|\leq c|\rho_1(1+\varepsilon)|^g,\quad i=1,\dots,n.
\end{equation*}
holds. This proves the first part of our theorem. If all roots $\rho_1,\dots,\rho_d$ with $\rho_1=\dots=\rho_d$ are simple, then
\begin{equation*}
|B_0^{(pg)}x_i-B_i^{(pg)}|\leq c|\rho_1|^g,\quad i=1,\dots,n.
\end{equation*}
The second part of the proof can be found in \cite[Theorem 59]{Sch2000}.
$\hfill\fbox{}$
\section*{Acknowledgements}
First, I thank my adviser Prof. Dr. Stefan M\"uller-Stach for picking out "multidimensional continued fractions" as topic of my diploma thesis and his helpful comments. Last but not least, I sincerely thank Prof. Dr. Fritz Schweiger for proofreading parts of this article and his valuable suggestions.
\bibliographystyle{alpha}
\bibliography{/Users/Kerli/Desktop/Selmer/SelmerEng}

\begin{thebibliography}{M{\"o}n54}

\bibitem[Bre81]{Bre1981}
Arne~Johan Brentjes.
\newblock {\em {Multi-Dimensional Continued Fractions}}.
\newblock Mathematisch Centrum, 1981.

\bibitem[M{\"o}n54]{Mon1954}
Rudolf M{\"o}nkemeyer.
\newblock {{\"U}ber Fareynetze in n Dimensionen}.
\newblock {\em Mathematische Nachrichten}, 11, 1954.

\bibitem[Pan08]{Pan2008}
Giovanni Panti.
\newblock {Multidimensional continued fractions and a Minkowski function}.
\newblock {\em Monatshefte f{\"u}r Mathematik}, (154):247--264, 2008.

\bibitem[Per07]{Per1907}
Oskar Perron.
\newblock {Grundlagen f{\"u}r eine Theorie des Jacobischen
  Kettenbruchalgorithmus}.
\newblock {\em Mathematische Annalen}, 64:1--76, 1907.

\bibitem[Sch00]{Sch2000}
Fritz Schweiger.
\newblock {\em {Multidimensional Continued Fractions}}.
\newblock Oxford University Press, 2000.

\bibitem[Sch05]{Sch2005}
Fritz Schweiger.
\newblock {Periodic Multiplicative Algorithms of Selmer Type}.
\newblock {\em Integers: Electronic Journal of Combinatorial Number Theory},
  5(1):1--9, 2005.

\end{thebibliography}
\end{document}